\documentclass[reqno]{amsart}

\usepackage[absolute,overlay]{textpos}
\usepackage{eso-pic}
\usepackage[a4paper,margin=1.25in]{geometry}
\usepackage{mathrsfs}
\usepackage{amsfonts,  amsmath}
\usepackage{graphicx}
\usepackage{enumerate,  enumitem,  amscd}
\usepackage{amsthm,  amssymb,  epsfig,  hyperref,  amsmath}  
\usepackage{mathtools}
\usepackage{tikz}
\usepackage{tikz-cd}
\usetikzlibrary{arrows,automata,arrows.meta}
\usepackage{quiver}
\usepackage{subcaption}
 
\tikzset{
  vx/.style  = {circle, draw=black, fill=blue!20, line width=0.6pt,
                minimum size=8pt, inner sep=0pt},
  ed/.style  = {draw=black!60, line width=0.5pt},
  lbl/.style = {font=\tiny},
  ttl/.style = {font=\small}
}

\newcommand{\pres}[3]{\textnormal{#1} \langle #2 \mid #3 \rangle}

\newcommand{\Z}{\mathbb{Z}}

\newcommand{\ab}{\mathrm{ab}}

\DeclareMathOperator{\BS}{BS}
\DeclareMathOperator{\BG}{BG}

\newtheorem*{theorem*}{Theorem} 
\numberwithin{theorem}{section}

\newtheorem*{corollary*}{Corollary}

\newtheorem*{proposition*}{Proposition}

\theoremstyle{definition}

\numberwithin{example}{section}
\newtheorem*{question*}{Question}
\numberwithin{question}{section}

\newtheorem*{remark*}{Remark}

\newcommand{\lastn}{0}
\InputIfFileExists{buildstamp.tex}{}{}

\newcounter{buildn}
\stepcounter{buildn}

\newcommand{\todaystamp}{\the\year.\ifnum\month<10 0\fi\the\month.\ifnum\day<10 0\fi\the\day}
\newcommand{\buildstamp}{\todaystamp\ \ v\thebuildn}

\newwrite\stampfile
\immediate\openout\stampfile=buildstamp.tex
\immediate\write\stampfile{\string}
\immediate\closeout\stampfile

\begin{document}

\title[A problem of Olshanskii]{The Baumslag--Gersten group and a problem of Olshanskii}
\author{Carl-Fredrik Nyberg-Brodda}
\address{June E Huh Center for Mathematical Challenges, Korea Institute for Advanced Study (KIAS), Seoul 02455, Korea}
\email{cfnb@kias.re.kr}

\thanks{The author is supported by KIAS Individual Grant HP094701 at Korea Institute for Advanced Study, and by the Mid-Career Researcher Program (RS-2023-00278510) through the National Research Foundation of Korea.}

\date{\today}


\begin{abstract} 
We prove that a certain representation of the Baumslag--Gersten one-relator group $\BG(1,2)$ by germs of continuous functions is not faithful. This gives a negative answer to a problem of A.~Yu.~Olshanskii from 2010 (Problem 17.99 in the Kourovka Notebook). 
\end{abstract}

\maketitle

\noindent Let $\mathcal{B} = \pres{}{a,b}{(bab^{-1})a(bab^{-1})^{-1} = a^2}$. This group was introduced by Baumslag \cite{Baumslag1969} in 1969, and is sometimes referred to as the Baumslag--Gersten group $\BG(1,2)$. It is an HNN extension of the Baumslag--Solitar group $\BS(1,2) = \pres{}{a,c}{cac^{-1}=a^2}$, with stable letter $b$ conjugating $\langle a \rangle \to \langle c \rangle$. A.\ Yu.\ Olshanskii noted that the functions $f(x) = 2x$ and $g(x) = 2^x$ satisfy the defining relation of $\mathcal{B}$, when composition is taken in the group of germs of monotonically increasing-to-$\infty$ continuous functions on $(0, \infty)$, i.e.\ where two functions are identified if they agree for sufficiently large arguments. Thus there is a representation $\mathcal{B} \to \langle f, g \rangle$ defined by $a \mapsto f$ and $b \mapsto g$. Olshanskii, in Problem 17.99 of the Kourovka Notebook \cite{Kourovka2022} (17th Issue, 2010), asked the natural question of whether this representation is faithful. We give a negative answer to this question by exhibiting a non-trivial element in the kernel of the representation. 

We directly construct the element. Consider the elements
\[
c = bab^{-1}, \quad d = c^{-1}ac, \quad s = bdb^{-1}, \quad \text{and} \quad t = sas^{-1}.
\]
In the above representation, with $g^{-1}$ represented by $x \mapsto \log_2(x)$ for large $x$, these map to
\[
C = 2^{2\log_2(x)} = x^2, \quad D = \sqrt{2x^2} = x\sqrt{2}, \quad S = 2^{\log_2(x)\sqrt{2}} = x^{\sqrt{2}}, \quad T = \left(2x^{\frac{1}{\sqrt{2}}}\right)^{\sqrt{2}} = 2^{\sqrt{2}}x
\]
respectively. In particular $[f, T]=1$, since $f(x) = 2x$ and $T(x) = 2^{\sqrt{2}}x$ are scalar multiples of $x$. However, as we now show, we have $[a, t] \neq 1$ in $\mathcal{B}$, proving that the representation is not faithful. 

Indeed, writing out the word $t$ in the letters $a, b, d$, we see that 
\begin{equation}\label{Eq:word}
[a,t] = ata^{-1}t^{-1} = a(bdb^{-1})a(bd^{-1}b^{-1})a^{-1}(bdb^{-1})a^{-1}(bd^{-1}b^{-1}).\tag{{$\star$}}
\end{equation}
Since the stable letter of the HNN extension above is $b$, we see that the only place that Britton pinches may appear is either as $b^{-1}a^{\pm 1}b$, or else as $bd^{\pm 1}b^{-1}$. The first type is not a pinch, since $a^{\pm 1} \not\in \langle c \rangle$ in $\BS(1,2)$. This may be seen by using the isomorphism $\BS(1,2) \cong \Z[\frac{1}{2}]\rtimes \Z$, where $n \in \Z$ acts via multiplication by $2^n$, and $a \mapsto (1,0)$ and $c \mapsto (0,1)$. For the second type, since $d = c^{-1}ac$ maps to $(\frac{1}{2},0)$, we have $d^{\pm 1}\not\in \langle a \rangle$, so  $bd^{\pm 1}b^{-1}$ is also not a pinch. Hence the word \eqref{Eq:word} contains no non-trivial pinches. By Britton's Lemma, we have $[a,t]\neq 1$ in $\mathcal{B}$. \qed

It would be interesting to understand the structure of the finitely generated group $\langle f, g \rangle$. In particular, is it finitely presented? Note that it is not hard to show that $\langle f, g \rangle^{\ab} \cong \Z$.

\bibliographystyle{amsalpha}
\bibliography{kourovka17-99.bib}

@book{Kourovka2022,
     title = {Unsolved Problems in Group Theory: The {K}ourovka {N}otebook},
    editor = {V. D. Mazurov and E. I. Khukhro},
    edition = {17th},
    year = {2010},
    address = {Novosibirsk},
    publisher = {Sobolev Institute of Mathematics}
}

@article {Baumslag1969,
    AUTHOR = {Baumslag, Gilbert},
     TITLE = {A non-cyclic one-relator group all of whose finite quotients
              are cyclic},
   JOURNAL = {J. Austral. Math. Soc.},
  FJOURNAL = {J. Austral. Math. Soc.},
    VOLUME = {10},
      YEAR = {1969},
     PAGES = {497--498},
}

\end{document}